\documentclass{amsart}
\usepackage{amssymb}
\usepackage{pdflscape}
\usepackage{amscd}
\usepackage{fancyhdr}
\usepackage{graphicx}

\usepackage[usenames,dvipsnames]{xcolor} 
\usepackage{tikz, tikz-3dplot, pgfplots}
\usepackage{tkz-graph}
\usepackage{tikz-cd}
\usetikzlibrary{positioning, calc, arrows, decorations.markings, decorations.pathreplacing}
\tikzset{
  symbol/.style={
    draw=none,
    every to/.append style={
      edge node={node [sloped, allow upside down, auto=false]{$#1$}}}
  }
}
\tikzcdset{scale cd/.style={every label/.append style={scale=#1},
    cells={nodes={scale=#1}}}}

\newtheorem{theorem}{Theorem}[section]
\newtheorem{corollary}[theorem]{Corollary}
\newtheorem{lemma}[theorem]{Lemma}
\newtheorem{proposition}[theorem]{Proposition}

\begin{document}
\title{Spherical representation spaces of quivers}
\author{Jennifer Müller, Markus Reineke}
\begin{abstract} We consider representation spaces of quivers, together with their base change action, and classify the spherical varieties among them.
\end{abstract}
\maketitle

\parindent0pt

\section{Introduction}

Spherical varieties, that is, varieties with a reductive group action for which a Borel subgroup acts with a dense orbit, form a particularly rewarding class of spaces, in particular since an extensive classification theory exists (see, for example, \cite{CF,Knop,KV,Lo,Lu,LV}). It is thus desirable to identify the spherical ones among natural classes of reductive group actions on varieties.\\[1ex]
A rather special, but particularly interesting, class of such actions is provided by the representation spaces of quivers with their base change action, for which the orbits naturally identify with isomorphism classes of quiver representations. The study of these actions can thus be enriched with representation-theoretic and homological techniques, leading to a fruitful theory of the geometry of orbits, invariants and quotients \cite{KacM,KR}.\\[1ex]
In the present work we classify the spherical varieties among the representation spaces of quivers. The main result, Theorem \ref{main}, states that these correspond to the cases where the quiver with its dimension vector is given as a connected sum, along entries $1$ of the dimension vector, of two special types of dimension vectors for quivers of Dynkin type $A_3$.\\[1ex]
The proof of this classification is quiver-theoretic, mainly relying on the theorem of Kac on indecomposable representations \cite{Kac}, and properties of infinite root systems.  In principle, our main result could also be obtained from the general classification of spherical linear representations of Kac, Brion, Benson, Radcliff, Leahy \cite[Section 5]{Knop}. However, the authors think that a concise quiver-theoretic derivation is  preferable, and the methods developed for this should be useful for a general study of Borel orbits in quiver representation spaces.\\[1ex]
The background on quiver representations, together with a key classical result of Brion and Vinberg on sphericity \cite{Brion,Vinberg}, are reviewed in Section \ref{pre}. Using a standard fibre bundle construction (Lemma \ref{bundle}) \cite{BR}, we identify in Section \ref{extend} Borel orbits in quiver representation spaces with certain orbits for the base change group in an extended quiver setting, for which finiteness can be characterized by a root-theoretic criterion (Proposition \ref{crit}). This allows us to classify the spherical representation spaces in Section \ref{classify}, by excluding certain quiver settings with infinitely many Borel orbits corresponding to imaginary roots of extended Dynkin diagrams.\\[1ex]
We note that Shmelkin \cite{S} has studied the much finer problem when a single orbit in a quiver representation space is spherical, and has developed a general approach to this problem. However, there seems to be no direct way to apply the results of \cite{S} to the present problem. Namely, identifying dense orbits in representation varieties (with respect to the base change action), and describing the decomposition behaviour of the corresponding representations (which is a requirement for the criteria of \cite{S}) leads to the intricate recursive combinatorics of infinite root systems studied in \cite{Scho}.\\[2ex]
{\bf Acknowledgments:} The authors would like to thank S.~Cupit-Foutou for references to the literature, F.~Knop for pointing out the classification in \cite{Knop}, and an anonymous referee for suggesting several improvements.

\section{Prerequisites}\label{pre}

We work over an algebraically closed base field $K$.  Let $Q$ be a finite quiver with set of vertices $Q_0$ and set of arrows $Q_1$, where arrows are written $\alpha:i\rightarrow j$ for $i,j\in Q_0$. Let ${\bf d}=(d_i)_i\in\mathbb{N}Q_0$ be a dimension vector for $Q$, and fix $K$-vector spaces $V_i$ of dimension $d_i$ for all $i\in Q_0$. The affine space
 $$R_{\bf d}(Q)=\bigoplus_{\alpha:i\rightarrow j}{\rm Hom}_K(V_i,V_j)$$
 encodes representations of $Q$ on the vector spaces $V_i$. The connected reductive algebraic group $$G_{\bf d}=\prod_{i\in Q_0}{\rm GL}(V_i)$$ acts on $R_{\bf d}(Q)$ via
 $$(g_i)_i\cdot(f_\alpha)_\alpha=(g_jf_\alpha g_i^{-1})_{\alpha:i\rightarrow j},$$
 and the orbits for this action correspond bijectively to the isomorphism classes of $K$-representations of $Q$ of dimension vector ${\bf d}$. The pair $(Q,{\bf d})$ is called a quiver setting.\\[1ex]
 For the remainder of this section, let $Q$ be acyclic. The Euler form of $Q$ is the (non-symmetric) bilinear form $\langle\_,\_\rangle_Q$ on $\mathbb{Z}Q_0$ given by
 $$\langle{\bf d},{\bf e}\rangle_Q=\sum_{i\in Q_0}d_ie_i-\sum_{\alpha:i\rightarrow j}d_ie_j.$$
 We denote its symmetrization by $(\_,\_)=(\_,\_)_Q$. The Weyl group $W_Q$ of $Q$ is the subgroup of ${\rm Aut}(\mathbb{Z}Q_0)$ generated by the reflections
 $$s_i({\bf d})={\bf d}-({\bf d},{\bf i})\cdot{\bf i}$$
 for $i\in Q_0$, where ${\bf i}\in\mathbb{Z}Q_0$ denotes the $i$-th coordinate vector. The union $\Delta^{\rm re}$ of $W_Q$-orbits of the ${\bf i}$ for $i\in Q_0$ is called the set of real roots. We denote by $F_Q$ the set of all $0\not={\bf d}\in\mathbb{N}Q_0$ which have connected support and fulfill $({\bf d},{\bf i})\leq 0$ for all $i\in Q_0$, called the {fundamental domain} of $Q$. Its $W_Q$-saturation $\Delta^{\rm im}=W_QF_Q$ is called the set of imaginary roots.  We denote by $\Delta^{{\rm re},+}$ (resp.~$\Delta^{{\rm im},+}$) the intersection of $\Delta^{\rm re}$ (resp.~$\Delta^{\rm im})$ with $\mathbb{N}Q_0$, called the set of positive real (resp.~imaginary) roots. For application in Section \ref{classify}, we note that every imaginary root ${\bf d}$ admits a root ${\bf e}\in F_Q$ in the same Weyl group orbit such that ${\bf e}\leq{\bf d}$. The theorem of Kac relates the infinite root system to the indecomposable representations of $Q$:
 
 \begin{theorem}\cite{Kac} There exists an indecomposable representation of $Q$ of dimension vector ${\bf d}$ if and only if ${\bf d}$ is a positive real or imaginary root. There exist infinitely many isomorphism classes of such representations if and only if ${\bf d}$ is imaginary.
 \end{theorem}
 
 This key result allows us to state a root-theoretic criterion for finiteness of orbits in representation spaces. We will not use the criterion itself, but derive a refinement (Proposition \ref{crit}) relevant for our setup in the next section.
 
 \begin{proposition}\label{precrit} The group $G_{\bf d}$ acts on $R_{\bf d}(Q)$ with finitely many orbits if and only if there is no positive imaginary root ${\bf e}\leq{\bf d}$.
 \end{proposition}
 
 \proof Assume that ${\bf e}\leq{\bf d}$ is a positive imaginary root. Then, by the previous theorem, there exists an infinite family $(U_\alpha)_\alpha$ of pairwise non-isomorphic indecomposable representations of $Q$ of dimension vector ${\bf e}$. Denoting by $S_{{\bf d}-{\bf e}}$ the representation of dimension vector ${\bf d}-{\bf e}$ in which all arrows are represented by zero maps, we thus find an infinite family $(U_\alpha\oplus S_{{\bf d}-{\bf e}})_\alpha$ of pairwise non-isomorphic representations of $Q$ of dimension vector ${\bf d}$, which correspond to infinitely many $G_{\bf d}$-orbits in $R_{\bf d}(Q)$, a contradiction.\\[1ex]
 Conversely, assume that no such ${\bf e}$ exists.
By the Krull-Remak-Schmidt theorem, every $V\in R_{\bf d}(Q)$ admits a decomposition $$V=U_1\oplus\ldots\oplus U_s$$ into indecomposables $U_k$, and the dimension vectors of the $U_k$ are real roots. Thus there are only finitely many isomorphism classes of such $V$. This proves finiteness of the number of $G_{\bf d}$-orbits in $R_{\bf d}(Q)$.\qed\\[1ex]
Our approach to sphericity of $R_{\bf d}(Q)$ will be based on the following result:

\begin{theorem}\label{finite}\cite{KV} Let $X$ be an irreducible algebraic variety with an action of a connected reductive algebraic group $G$, and let $B$ be a Borel subgroup of $G$. Then $X$ admits an open $B$-orbit if and only if $B$ acts on $X$ with finitely many orbits.
\end{theorem}

\section{Leg-extended quivers}\label{extend}

 We fix complete flags $F_i^*$ in $V_i$ for all $i\in Q_0$, and denote by $B_{\bf d}\subset G_{\bf d}$ the Borel subgroup fixing all $F_i^*$.  We are interested in a description of the $B_{\bf d}$-orbits in $R_{\bf d}(Q)$.\\[1ex]
 We define the {\it leg-extended quiver} ${Q}_{\bf d}$ by
 \begin{itemize}
 \item vertices $(i,k)$ for $i\in Q_0$ and $k=1,\ldots,d_i$,
 \item arrows $\beta_{i,k}:(i,k)\rightarrow (i,k+1)$ for $i\in Q_0$ and $k=1,\ldots,d_i-1$,
 \item arrows $\alpha:(i,d_i)\rightarrow(j,d_j)$ for $\alpha:i\rightarrow j$ an arrow in $Q$.
 \end{itemize}
 
 We also define $\overline{Q}_{\bf d}$ by omitting the arrows $\alpha$ in $Q_{\bf d}$. We define a dimension vector $\widehat{\bf d}$ for ${Q}_{\bf d}$ by $$\widehat{d}_{(i,k)}=k$$ for $i\in Q_0$ and $k=1,\ldots,d_i$. We call a representation $V$ of ${Q}_{\bf d}$ of dimension vector $\widehat{\bf d}$ {\it of flag type} if all $\beta_{i,k}$ are represented by injective maps. In this case, up to isomorphism, we can assume that $V_{(i,k)}=F_i^*$ for all $i$ and $k$, and that $V_{\beta_{(i,k)}}$  equals the inclusion. We thus find a canonical bijection between the $G_{\widehat{\bf d}}$-orbits in the open subset $R_{\widehat{\bf d}}^{\rm flag}(Q_{\bf d})$ of $R_{\widehat{\bf d}}(Q_{\bf d})$ of representations of flag type and the $B_{\bf d}$-orbits in $R_{\bf d}(Q)$. More precisely, the restriction map $$p:R_{\widehat{\bf d}}^{\rm flag}(Q_{\bf d})\rightarrow R_{\widehat{\bf d}}(\overline{Q}_{\bf d})$$ is $G_{\widehat{\bf d}}$-equivariant onto a homogeneous space with stabilizer $B_{\bf d}$, with fibre isomorphic to $R_{\bf d}(Q)$. We conclude:
 
 \begin{lemma}\label{bundle} We have a $G_{\widehat{\bf d}}$-equivariant isomorphism
 $$R_{\widehat{\bf d}}^{\rm flag}(Q_{\bf d})\simeq G_{\widehat{\bf d}}\times^{B_{\bf d}}R_{\bf d}(Q).$$
 In particular, $B_{\bf d}$-orbits in $R_{\bf d}(Q)$ correspond bijectively to $G_{\widehat{\bf d}}$-orbits in $R_{\widehat{\bf d}}^{\rm flag}(Q_{\bf d})$.
 \end{lemma}
 
We will now refine Proposition \ref{precrit} to the case of representations of flag type for a leg-extended quiver. For this, we call a dimension vector ${\bf e}\leq\widehat{\bf d}$ for $Q_{\bf d}$ {\it of flag type} if $$e_{(i,1)}\leq\ldots\leq e_{(i,d_i)}$$ for all $i\in Q_0$. We call ${\bf e}$ {\it gentle} if both ${\bf e}$ and $\widehat{\bf d}-{\bf e}$ are of flag type. Equivalently, ${\bf e}$ is gentle if and only if $$e_{(i,k+1)}-e_{(i,k)}\in\{0,1\}$$ for all $i\in Q_0$ and $k<d_i$. 

\begin{lemma} If ${\bf e}$ is an imaginary root for $Q_{\bf d}$ and $V$ is an indecomposable representation of dimension vector ${\bf e}$, then $V$ is of flag type.
\end{lemma}

\proof Let ${\bf e}$ be an imaginary root for $Q_{\bf d}$, and let $V$ be an indecomposable representation of dimension vector ${\bf e}$. If $V$ is not of flag type, we can find $i\in Q_0$ such that $V_{\beta_{(i,k)}}$ is not injective for some $k<d_i$; let $k$ be minimal with this property. Choose a non-zero vector $v_k$ such that $$V_{\beta_{(i,k)}}(v_k)=0.$$ Choose $l\leq k$ minimal with the property that $$v_k=V_{\beta_{(i,k-1)}}\ldots V_{\beta_{(i,l)}}v_l,$$ and define $$v_p=V_{\beta_{(i,p-1)}}\ldots V_{\beta_{(i,l)}}v_l$$ for $l<p<k$. The $v_p$ for $l\leq p\leq k$ clearly span a non-zero subrepresentation $U$ of $V$. We claim that this is even a direct summand. Namely, we first define $$W_{(i,p)}=V_{(i,p)}$$ for all $p<l$. For $l\leq p\leq k$, inductively, we choose a complement $W_{(i,p)}$ to $v_p$ containing $V_{\beta_{(i,p-1)}}W_{(i,p-1)}$, which is possible by the minimality properties of $k$ and $l$ and by the choice of the $v_p$. We define $W_x=V_x$ for all other vertices $x$ of $Q_{\bf d}$. Then $W$ defines a subrepresentation of $V$ complementing $U$, thus $V=U$ by indecomposability. But the dimension vector of $U$ is clearly a real root, yielding a contradiction.\qed\\[1ex]
The proof of the following follows closely the proof of Proposition \ref{precrit}.
 
 \begin{proposition}\label{crit} For acyclic $Q$, the group $G_{\widehat{\bf d}}$ acts on $R_{\widehat{\bf d}}^{\rm flag}(Q)$ with finitely many orbits if and only if there is no gentle imaginary root ${\bf e}\leq{\widehat{\bf d}}$.
 \end{proposition}
 
 \proof Assume that ${\bf e}\leq{\widehat{\bf d}}$ is a gentle imaginary root. Then, there exists an infinite family $(U_\alpha)_\alpha$ of pairwise non-isomorphic indecomposable representations of $Q$ of dimension vector ${\bf e}$, which are of flag type by the previous proposition. Since $\widehat{\bf d}-{\bf e}$ is again of flag type, we can choose a representation $S$ of flag type of this dimension vector. We thus find an infinite family $(U_\alpha\oplus S)_\alpha$ of pairwise non-isomorphic representations of flag type of $Q_{\bf d}$ of dimension vector $\widehat{\bf d}$, which correspond to infinitely many $G_{\widehat{\bf d}}$-orbits in $R_{\widehat{\bf d}}(Q_{\bf d})$, a contradiction.\\[1ex]
 Conversely, assume that no such ${\bf e}$ exists, but that there are infinitely many isomorphism classes of representations of flag type.
By the Krull-Remak-Schmidt theorem, every $V\in R_{\widehat{\bf d}}^{\rm flag}(Q_{\bf d})$ admits a decomposition $$V=U_1\oplus\ldots\oplus U_s$$ into indecomposables $U_k$ which are again of flag type. If the dimension vector ${\bf e}$ of one of them, say $U_1$, is an imaginary root, it is gentle, since both $U_1$ and $U_2\oplus\ldots\oplus U_s$ are of flag type, a contradiction. Thus, the dimension vector of any such $U_k$ is a real root, and there are only finitely many isomorphism classes of such $V$. This proves finiteness of the number of $G_{\widehat{\bf d}}$-orbits in $R_{\widehat{\bf d}}^{\rm flag}(Q)$.\qed

\begin{corollary} For ${\bf d}\in\mathbb{N}Q_0$, the representation space $R_{\bf d}(Q)$ is spherical if and only if there is no gentle imaginary root ${\bf e}\leq\widehat{\bf d}$ for $Q_{\bf d}$.
\end{corollary}

\proof Sphericity of $R_{\bf d}(Q)$ is equivalent to $B_{\bf d}$ acting with finitely many orbits on $R_{\bf d}$ by Theorem \ref{finite}. This is equivalent to $G_{\widehat{\bf d}}$ acting with finitely many orbits on $R_{\widehat{\bf d}}^{\rm flag}(Q_{\bf d})$, and the previous proposition shows the claim.\qed

\section{Classification}\label{classify}

To approach the classification of spherical representation varieties, we first show that we can restrict to excluding minimal non-spherical situations.

\begin{lemma} If ${\bf d}\leq{\bf e}$ are dimension vectors for $Q$ and $R_{\bf d}(Q)$ is not spherical, then $R_{\bf e}(Q)$ is not spherical.
\end{lemma}

\proof We have an embedding of quivers $$\iota:Q_{\bf d}\rightarrow Q_{\bf e},\;\;\; (i,k)\mapsto(i,k+e_i-d_i)$$ for $i\in Q_0$ and $k\leq d_i$. It induces an embedding of Weyl groups $$\iota:W_{Q_{\bf d}}\rightarrow W_{Q_{\bf e}}$$ and a $W_{Q_{\bf }}$-equivariant embedding $$\iota:\mathbb{Z}(Q_{\bf d})_0\rightarrow\mathbb{Z}(Q_{\bf e})_0.$$ It is easily verified that $\iota$ maps $F_{Q_{\bf d}}$ to $F_{Q_{\bf e}}$, and thus imaginary roots for $Q_{\bf d}$ to imaginary roots for $Q_{\bf e}$. Obviously, $\iota$ also preserves gentleness. The claim then follows from the previous corollary.\qed\\[1ex]
We say that a quiver $Q$ is a {\it connected sum} of two full subquivers $Q'$ and $Q''$ along a vertex $i\in Q_0$ if $$(Q')_0\cup (Q'')_0=Q_0\mbox{ and }(Q')_0\cap(Q'')_0=\{i\}\mbox{ and }(Q')_1\cap (Q'')_1=\emptyset.$$
In this case, we write $$Q=Q'\#_iQ''.$$ Given a dimension vector ${\bf d}$ for $Q$, we write
$$(Q,{\bf d})=(Q',{\bf d}|_{Q'})\#_i(Q'',{\bf d}|_{Q''}).$$
We call such a connected sum of quiver settings {\it thin} if $d_i=1$. In the following, the symbol $\leftrightarrow$ denotes an arrow of either orientation.

\begin{theorem}\label{main}
For a connected quiver setting $(Q,\boldsymbol{d})$, the representation space $R_{\boldsymbol{d}}(Q)$ is a spherical variety under the action of $B_{\boldsymbol{d}}\subset G_{\boldsymbol{d}}$ if and only if $Q$ is of Dynkin type $A_1$ or $A_2$ or if $(Q,\boldsymbol{d})$ is a thin connected sum of quiver settings with $Q$ of the form 
$\begin{tikzcd}[column sep=small]
	1 \arrow[r,leftrightarrow] & 2 \arrow[r,leftrightarrow] & 3
\end{tikzcd}$,
and $\boldsymbol{d}=(m,n,1)$ or $\boldsymbol{d}=(m,2,n)$.

\end{theorem}


\proof First assume $R_{\bf d}(Q)$ is spherical. We reduce to the claimed situation in several steps, always using the previous lemma.
\begin{enumerate}
\item If $Q$ contains a cycle, then $R_{\bf d}(Q)$ even admits infinitely many $G_{\bf d}$-orbits, since there are already infinitely many isomorphism classes of representations of dimension vector $(1,\ldots,1)$ for any cycle quiver. Thus we can assume $Q$ to be a tree.
\item With the same argument applied to the Kronecker quiver, we can assume $Q$ to have no parallel arrows between vertices.
\item Assume that $Q$ cannot be written as a proper thin connected sum.
\item Then there is no vertex $i$ incident with more than two other vertices: otherwise, we can assume $d_i\geq 2$ and $i$ being connected to $j_1, j_2 ,j_3$, and the gentle imaginary root $${\bf (i,1)}+2{\bf (i,2)}+{\bf (j_1,1)}+{\bf (j_2,1)}+{\bf (j_3,1)}$$ for a $\widetilde{D}_4$-quiver contradicts sphericity. Thus we can assume $Q$ to be of type $A_n$.
\item If $n\geq 4$, we can find a subquiver $i_1\leftrightarrow i_2\leftrightarrow i_3\leftrightarrow i_4$ such that $d_{i_2},d_{i_3}\geq 2$, and the gentle imaginary root $${\bf (i_1,1)}+{\bf (i_2,1)}+2{\bf (i_2,2)}+{\bf (i_3,1)}+2{\bf (i_3,2)}+{\bf (i_4,1)}$$ for a $\widetilde{D}_5$ quiver contradicts sphericity. Thus $Q$ is at most of type $A_3$.
\item Say $Q$ is of the form $i_1\leftrightarrow i_2\leftrightarrow i_3$. The gentle imaginary root $\widehat{(2,3,2)}$ for a $\widetilde{E}_6$-quiver shows that either $d_{i_2}\leq 2$, or $d_{i_1}\leq 1$, or $d_{i_3}\leq 1$.
\end{enumerate}
This finishes the first part of the proof.\\[1ex]
To prove the converse direction, we first show that the $A_3$-situations of the theorem are indeed spherical. We exemplify this for the case $\boldsymbol{d}=(m,n,1)$ by using Corollary 3.4 and direct calculations.
       \begin{figure}[h!]
        \centering
        \[\begin{tikzcd}[sep=small]
            e_{(1,d_1)} \arrow[r,leftrightarrow] & e_{(2,d_2)} \arrow[r,leftrightarrow] & e_{(3,1)} \\
            e_{(1,d_1-1)} \arrow[u,hookrightarrow] & e_{(2,d_2-1)} \arrow[u,hookrightarrow] & \\
            \vdots \arrow[u,hookrightarrow] & \vdots \arrow[u,hookrightarrow] & \\
            e_{(1,1)} \arrow[u,hookrightarrow] & e_{(2,1)} \arrow[u,hookrightarrow] & 
        \end{tikzcd}\]
        \caption{Labelling of the vertices in the case $\boldsymbol{d}=(m,n,1)$}
        \label{fig:placeholder}
    \end{figure} \\
    So assume $\boldsymbol{e}\leq \boldsymbol{\widehat{d}}$ is an imaginary root for $Q_{\boldsymbol{d}}$. We can assume $e_{(3,1)}=1$, otherwise the support of ${\bf e}$ is of type $A$, and $\boldsymbol{e}$ is a real root. We may also assume $m$ and $n$ to be minimal in the sense that $e_{(1,1)}$ and $e_{(2,1)}$ are non-zero, thus $e_{(1,1)}=e_{(2,1)}=1$. Moreover, by applying reflections, we may assume that $\boldsymbol{e}\in F_{Q_{\bf d}}$. This allows us to use the conditions $(\boldsymbol{e},\boldsymbol{(i,k)}) \leq 0$ for all $\boldsymbol{(i,k)} \in Q_{\boldsymbol{d}}$. We first find
    \[(\boldsymbol{e},\boldsymbol{(i,1)}) = 2e_{(i,1)} - e_{(i,2)} \leq 0 \;\Leftrightarrow \;e_{(i,2)} \geq 2\]
    for $i\in \{1,2\}$. But we also have $e_{(i,2)}\leq d_{(i,2)}=2$, thus $e_{(i,2)} = 2$. By a quick induction, we find $e_{(i,k)}=k$ for all $k\leq d_i$. Additionally, we find
    \[(\boldsymbol{e},\boldsymbol{(1,d_1)}) = 2 e_{(1,d_1)} - e_{(1,d_1-1)} - e_{(2,d_2)} = 2d_1-d_1+1-d_2 \leq 0,\]
    thus $d_1+1 \leq d_2$, and
    \[(\boldsymbol{e},\boldsymbol{(2,d_2)}) = 2 e_{(2,d_2)} - e_{(2,d_2-1)} - e_{(1,d_1)} - e_{(3,1)} = 2d_2-d_2+1-d_1-1 \leq 0, \]
  thus $d_2 \leq d_1$, a contradiction.  The cases $\boldsymbol{d}=(m,n)$ and $\boldsymbol{d}=(m,2,n)$ are proven similarly.

Now we prove that sphericity is preserved under thin connected sums. First note that a connected sum of trees is again a tree by disjointness of the respective sets of arrows. So let $(Q,\boldsymbol{f})=(Q',\boldsymbol{c}) \#_i (Q'',\boldsymbol{d})$ be the thin connected sum along a vertex $i\in Q_0$ such that $R_{\boldsymbol{c}}(Q')$ is spherical, $Q''$ is of type $A_3$, and $\boldsymbol{d}=(m,n,1)$.
   Let $\boldsymbol{e} \leq \boldsymbol{\widehat{f}}$ be a positive root, which again we can assume to belong to $F_{Q_{\bf f}}$. We can then repeat the previous arguments word by word to arrive at the same contradiction, finishing the proof.\qed

\bibliography{Spherical}{}

\begin{thebibliography}{10}

\bibitem{BR}
Maria Bertozzi and Markus Reineke.
\newblock Momentum map images of representation spaces of quivers.
\newblock {\em J. Lie Theory}, 32(3):797--812, 2022.

\bibitem{Brion}
Michel Brion.
\newblock Quelques propri\'et\'es des espaces homog\`enes sph\'eriques.
\newblock {\em Manuscripta Math.}, 55(2):191--198, 1986.

\bibitem{CF}
St\'ephanie Cupit-Foutou.
\newblock Spherical varieties and perspectives in representation theory.
\newblock In {\em Representation theory---current trends and perspectives}, EMS
  Ser. Congr. Rep., pages 47--57. Eur. Math. Soc., Z\"urich, 2017.

\bibitem{Kac}
V.~G. Kac.
\newblock Infinite root systems, representations of graphs and invariant
  theory.
\newblock {\em Invent. Math.}, 56(1):57--92, 1980.

\bibitem{KacM}
Victor~G. Kac.
\newblock Root systems, representations of quivers and invariant theory.
\newblock In {\em Invariant theory ({M}ontecatini, 1982)}, volume 996 of {\em
  Lecture Notes in Math.}, pages 74--108. Springer, Berlin, 1983.

\bibitem{Knop}
Friedrich Knop.
\newblock Some remarks on multiplicity free spaces.
\newblock In {\em Representation theories and algebraic geometry ({M}ontreal,
  {PQ}, 1997)}, volume 514 of {\em NATO Adv. Sci. Inst. Ser. C: Math. Phys.
  Sci.}, pages 301--317. Kluwer Acad. Publ., Dordrecht, 1998.

\bibitem{KV}
Friedrich Knop and Bart Van~Steirteghem.
\newblock Classification of smooth affine spherical varieties.
\newblock {\em Transform. Groups}, 11(3):495--516, 2006.

\bibitem{KR}
H.~Kraft and Ch. Riedtmann.
\newblock Geometry of representations of quivers.
\newblock In {\em Representations of algebras ({D}urham, 1985)}, volume 116 of
  {\em London Math. Soc. Lecture Note Ser.}, pages 109--145. Cambridge Univ.
  Press, Cambridge, 1986.

\bibitem{Lo}
Ivan~V. Losev.
\newblock Proof of the {K}nop conjecture.
\newblock {\em Ann. Inst. Fourier (Grenoble)}, 59(3):1105--1134, 2009.

\bibitem{Lu}
D.~Luna.
\newblock Vari\'et\'es sph\'eriques de type {$A$}.
\newblock {\em Publ. Math. Inst. Hautes \'Etudes Sci.}, (94):161--226, 2001.

\bibitem{LV}
D.~Luna and Th. Vust.
\newblock Plongements d'espaces homog\`enes.
\newblock {\em Comment. Math. Helv.}, 58(2):186--245, 1983.

\bibitem{Scho}
Aidan Schofield.
\newblock General representations of quivers.
\newblock {\em Proc. London Math. Soc. (3)}, 65(1):46--64, 1992.

\bibitem{S}
D.~A. Shmelkin.
\newblock On spherical representations of quivers and generalized complexes.
\newblock {\em Transform. Groups}, 7(1):87--106, 2002.

\bibitem{Vinberg}
\`E.\~B. Vinberg.
\newblock Complexity of actions of reductive groups.
\newblock {\em Funktsional. Anal. i Prilozhen.}, 20(1):1--13, 96, 1986.

\end{thebibliography}
\bibliographystyle{plain}

\end{document}